\documentclass{amsart}

\usepackage[T1]{fontenc}
\usepackage[utf8]{inputenc}
\usepackage[main=english,russian]{babel}
\newcommand{\rus}[1]{\foreignlanguage{russian}{#1}}

\usepackage{graphicx}        

\usepackage[bottom]{footmisc}
\usepackage{url}

\usepackage{newtxtext}       %
\usepackage{newtxmath}       

\begin{document}

\title[Mathematics discovered, invented, and inherited]{A view from lockdown:\\  mathematics discovered, invented, and inherited}

\author[A.V.Borovik]{Alexandre Borovik}
\address{Department of Mathematics, University of Manchester, UK}
\email{alexandre$\gg$at$\ll$borovik.net}

\thanks{\copyright\ 2020 Alexadre Borovik under exclusive license to Springer Nature Switzerland AG 2020\\
for a  shorter version of the paper (roughly Sections $1$ to $4$) which appeared as:
Borovik A. (2020) A View from Lockdown: Mathematics Discovered, Invented, and Inherited. In: . Mathematics Online First Collections. Springer, Cham.  \url{https://doi.org/10.1007/16618_2020_6}.}

\date{28 December 2020}

\maketitle


\noindent
\textbf{\textsc{Abstract.}} The classical platonist / formalist dilemma in philosophy of mathematics can be expressed in lay
terms as a deceptively naive question:
\begin{center}
\emph{Is new mathematics discovered or invented?}
\end{center}
Using examples from my own mathematical work during the Coronavirus lockdown, I argue that there is also a third way:
new mathematics can also be\emph{ inherited}. And entering into possession, making it your own,  could be great fun.\\

\noindent
\textbf{\textsc{Executive Summary.}} This paper is written in lockdown and can serve as a testimony in support of the apparently  self-evident, but largely ignored  principle:

\begin{quote} The most important resource for a (pure) mathematician's research is uninterrupted time for thinking.
\end{quote}
University administrators and research funding bodies systematically ignore it, and the bureaucratic burden imposed by them strangulates mathematics research. A short breathing space provided by lockdown could make miracles.

\section{July 2020. Your best friend, the subconscious}

\begin{flushright}
\emph{The soul is silent. If it speaks at all it speaks in dreams.}\\
Louise Gl\"{u}ck
\end{flushright}

I confess, with some embarrassment, that my life in lockdown is comfortable and happy. I wake up at sunrise and take an hour long walk in the local park (conveniently, a wilderisation project), meeting on my way only foxes and birds -- among them the resident grey heron, \emph{Ardea cinerea}, an elegant and dignified bird.\footnote{A few days ago (that is, end of July) I have seen a relatively rare atmospheric phenomenon: a full arch double rainbow, of very intense colour, at the very moment when the rising sun was crossing the horizon. Some would perhaps see that as a good omen and symbol of hope, or a tribute to the National Health Service (badly painted rainbows are everywhere all over the country), but I instead started to construct, in my head, a geometric explanation of the old conundrum: why, in the two arches of a double rainbow, colours change in opposite orders: from blue to red in the inner, and from red to blue in the outer arch? I leave this problem to the readers as an exercise.}   After a light breakfast and coffee, I start doing mathematics, that is, I sit at my desk and look out of the window. This is a hard job, and I soon become tired, move to a sofa and take a nap. On waking up, I am refreshed, and return to mathematics -- and more often than not I have some new ideas for my work;  they came to me during my sleep. This cycle is repeated, with breaks for meals and tea. My wife Anna uses meals for briefing me about \textsc{Covid} and other news, I myself do not follow current affairs.

Perhaps at this point I have to touch on one of the best kept secrets of mathematics:
\begin{center}
\emph{mathematics is done in the subconscious.}
\end{center}
A mathematician has to maintain good relations with his or her subconscious. The subconscious is not a properly domesticated beast, but it responds well to attention and kindness. It is like our rabbit, Cadbury the Netherland Dwarf (one of the wilder breeds of pet rabbits). When he is in good spirits, Cadbury grooms me, combing with his incisors the skin on my arm, apparently trying to relieve me of my (non-existent, I hope) fleas -- this is a natural social behaviour of rabbits. While I doze, my subconscious combs the deepest recesses of  my memory for morsels of mathematics which could be relevant to, or just somehow associated with the mathematics that I am trying to do in my conscious state. The subconscious is a wordless creature and brings its catch to the surface as a kind of uncertain, instantly disappearing visual image akin to a single frame inserted in a film reel (the ``inverse vision'' as described by William Thurston \cite{Thurstion1994}). Then another miracle happens: someone or something else in my mind looks at the catch and says: ``well, this is \dots'' -- and gives the name, usually an already well-known term of mathematical language.\footnote{You may find more on that in my paper \cite[Section 6.1]{borovik-makers-users}.}

\indent A few days ago, a crucial ingredient of a proof on which I was working was brought that way to the surface after hibernating in the depths for four decades -- completely forgotten and never touched by me. Now, when I write this story, I am able to recall the particular paper where I first encountered this concept, and the monograph which I consulted to learn more about it.\footnote{For the mathematician reader: the paper and the monograph were the seminal paper by Hall and Higman \cite{hall-higman} and the classical book by Curtis and Reiner \cite{curtis-reiner}. And the concept was the \emph{enveloping algebra} of a representation. Hall and Higman revolutionised the finite group theory by observing that if $L/K$ and $M/L$ are sections in a finite group $G$ for $ 1 \leqslant K \lhd L \lhd M \leqslant   G$, with $K \lhd M$ and $L/K$ being an elementary abelian group of order $p^l$ for prime $p$, then the action of $M/L$ on $L/K$ by conjugation is a representation $M/L \longrightarrow {\rm GL}_l(\mathbf{F}_p)$  can be usefully studied by methods of representation theory. The images of elements from $M/L$ generate a subring (called the \emph{enveloping algebra}) in the matrix algebra ${\rm M}_{l\times l}(\mathbf{F}_{p})$. If the representation is irreducible, this subring is the matrix algebra ${\rm M}_{m\times m}(\mathbf{F}_{p^n})$ for $mn= l$ \cite[Lemma 70.5]{curtis-reiner}.\\[.5ex]
I made the following basic observation only now, in lockdown: the group $G$ is much more complex than it looks at the first glance because the general linear group ${\rm GL}_{m\times m}(\mathbf{F}_{p^n}) \subset {\rm M}_{m\times m}(\mathbf{F}_{p^n})$ lives naturally  and \emph{acts} inside $G$ (in terminology of model theory, it is \emph{interpretable} in $G$). Moreover, the groups ${\rm GL}_{m\times m}(\mathbf{F}_{p^n})$ are some of the better understood finite groups. This brought dramatic simplifications  into the model theoretic problem  I was working on.} All that had happened in about 1976, when I was an undergraduate student. I never touched the stuff since then.

I reached this state of nirvana only at the end of May, when the heroic attempts to teach online had been paused for summer. In lockdown before that, I lived comfortably, taught online quite productively, wrote papers, but did not achieve the wholeness of being that I experience now.

\section{Mathematics inherited}

And now I turn to the issue indicated in the title of my notes:
\begin{center} \emph{mathematics discovered, invented},\footnote{This philosophical dilemma -- discovered vs invented -- has interesting practical implications: a mathematical formula can be patented in the USA, but not in the UK. For American lawmakers and lawyers, the formula is invented, for British -- discovered.} \emph{and inherited}.
\end{center}
My story is about mathematics which is neither discovered nor invented, but inherited. Indeed I argue that
\begin{center}
\emph{new mathematics can also be inherited.}
\end{center}
Some years ago I made  my paper \cite{borovik-inherited} public where I described an example from my own work,  a convoluted pre-history of a few simple but powerful mathematical ideas and the way they were inherited, or ignored, or re-discovered by new generations of mathematicians. That old paper was too technical, but I borrow from it a couple of softer passages.

My lockdown episode fits into the same pattern of inheritance: I recovered from my memory a few elementary mathematical concepts which I learned decades ago and have not used for 40+ years; they belonged to canonical stuff, occasionally taught in senior undergraduate courses, found in textbooks and therefore securely fossilised. This happened because I worked on sufficiently hard problems in model-theoretic algebra (one of them, \cite[Problem B.38, p. 365]{borovik-nesin}, could be traced back to 1994) . I realised that I needed to get some good understanding of how the traditional old stuff could be used in the new non-traditional environment.

And I suddenly got a rather frivolous idea. I decided to experiment with the principle that I formulated  years ago, perhaps in my undergraduate years:

\begin{center}
\emph{mathematics can be done with a matchstick on a moldy wall of a prison cell.}
\end{center}

For many years, I maintained a list of rainy day problems, something that I could do without access to the literature -- I thought the list could be useful if I was confined to a bed in a hospital with just paper and a pencil to save me from boredom. I included one of these problems in my book \cite[Section 11.2]{borovik-microscope}: it is about development of Euclidean geometry from an alternative system of axioms -- something that, I had a good reason to believe, could indeed be done with a matchstick on a prison cell wall. Now, in lockdown, I decided to turn my research project into a matchstick exercise -- work on it  from the first principles,  without looking into any books or any external sources of information, and even not making notes on paper.

To my joy, I almost instantly started to feel that I was developing a much deeper insight into the problem than I would otherwise have had, and that the emerging proof was, as mathematicians would say, the `right' proof. It is now a cute little paper \cite{borovik-finite} (its wording is much more formal than my original sketch). It was incredible fun -- I laughed watching transformation of ideas.

\section{ Ernst Haeckel: Ontogeny recapitulates phylogeny}

My lockdown experience was an illustration of Haeckel'principle as expressed in the title of this section.

In the past
I was privileged to work with Israel Gelfand, one of the great mathematicians of the 20th century. He made a clear distinction between the two modes of work in mathematics expressed by Russian words `\emph{pr\textbf{I}dumyvanie}' and `\emph{pr\textbf{O}dumyvani}e', very similar and almost homophonic. The former means `inventing,' the latter `properly thinking through' and was used by Gelfand with the meaning
\begin{center}
\emph{`thinking through starting from the origins, fundamentals, first principles.'}
\end{center}
Gelfand valued \emph{pr\textbf{O}dumyvanie} more than \emph{pr\textbf{I}dumyvanie}, he was convinced that \emph{pr\textbf{O}dumyvanie} yielded deeper results.

I was lucky that I  followed Gelfand's advice and restricted myself to \emph{pr\textbf{O}dumyvanie} -- this gave me a few days of happiness and an immense intellectual joy.

\section{Inoculation against lockdown blues}

Finally, I have to explain what allows me to be happy under lockdown.

Mathematics is a proselytising cult, and by the age of 16 I swallowed its dogmata hook, line, and sinker, and became an unwavering convert.  Myself and  friends at FMSh, the Physics and Mathematics Boarding  School at the Novosibirsk University, knew that we were to become professional researchers in physics or mathematics, and, moreover, we knew that we had no other choice  because this was the only way available to us to maintain some degree of intellectual freedom. At that time, in Soviet specialist schools like FMSh this was a commonplace sentiment -- see \cite{Gerovitch2019}.

During my first year at university another colour was added to this vision of the world: the explicit understanding (shared by my friends) that
\begin{center}
\emph{mathematics was the best  escape route from reality.}
\end{center}

At the next stage of my professional development, when I was a kind of a postdoc, I got hold of the Russian translation of the novel \emph{Theophilus North} by Thornton Wilder \cite{Wilder1973}. I already knew earlier works by Wilder, and opened the book with some anticipation.

Theophilus, the narrator and protagonist of the novel, tells about himself at the very beginning of the novel:
\begin{quote}
At various times I had been afire with \textsc{Nine Life Ambitions} -- not necessarily successive, sometimes concurrent, sometimes dropped and later revived, sometimes very lively but under a different form and only recognized, with astonishment, after the events which had invoked them from the submerged depths of consciousness.
\end{quote}
He gives a curious list:
\begin{quote}
a saint, an anthropologist among primitive peoples, the archaeologist, the detective, the actor,  the magician, the lover, the rascal, -- and a free man.
\end{quote}
I immediately put the book aside and started my own list. The translator did not use the word `ambitions' (perhaps because Soviet people were not supposed to have ambitions, this word had negative connotations); the one used could be translated back as `role' or `field of activity'. This made my spontaneously produced list a bit more precise. Most entries are irrelevant now, but the last one matters: \emph{political exile}. Still a young man, I suddenly discovered that I was  prepared to face this fate.

Should I add anything else to claim that growing up in certain political environments is the best inoculation against lockdown blues later in life?

\section{August-October 2020. Nothing to report}

\begin{flushright}
\emph{In the spring of 1926 I resigned from my job.}\\
\emph{The first days following such a decision are like}\\
\emph{the release from a hospital after a protracted illness}.\\
\emph{One slowly learns how to walk again;}\\
\emph{slowly and wonderingly one raises one's head}.\\[1ex]
Thornton Wilder,\\
the opening lines of \emph{Theophilus North}.
\end{flushright}

I spent August and September drowning in the bureaucratic nightmare of retiring from my university, and October in a slow recovery. No Nirvana. Nothing to report.

However, this led me to formulating the \emph{Executive Summary} (please see the first page of this text). Also, please check my Disclaimer (at the end of the text, just before \emph{References}).

\section{November 2020. Memory and knowledge: taxonomy}

\begin{flushright}
\emph{We look at the world once, in childhood.}\\  \emph{The rest is memory.}\\
Louise Gl\"{u}ck
\end{flushright}

In the new lockdown, I returned to my mathematical work with the same respect to intuition and to the subconscious, and quite successfully, which prompted my decision to return to this my paper where I document my reflections on my mathematical work.

The story that I have told so far raises questions about the role of memory, and of the specific kind of memory: ability to trace the stuff that you know back  to the origins,  to the first instances when you learned it. Moreover,  tracing your knowledge  back to its origins inevitably  involves maintaining in the head some, maybe not well organised, but still a kind of a system for recording links with separate pieces of your knowledge.  The first readers of  my paper promptly asked  me how this was possible.

These aspects of human's memory are critically important for mechanisms of  ``inheritance'' in mathematics -- and for doing serious mathematics --  but, to the best of my knowledge,  are not discussed anywhere. I have to resort to my own experience of developing them and also to refer to the authority in applied philosophy -- Donald Rumsfeld. I heard, and was prepared to support, proposals to award to Rumsfeld the Ig Nobel prize\footnote{About the Ig Nobel Prizes: \texttt{https://www.improbable.com/ig-about/}.} in epistemology\footnote{The only Ig Nobel prize for work of philosophical nature was awarded in 2011 (under the Literature category) to  John Perry of Stanford University for his \emph{Theory of Structured Procrastination}, which states: ``To be a high achiever, always work on something important, using it as a way to avoid doing something that's even more important'' \cite{Perry1996,Perry2012}; see also \texttt{http://www.structuredprocrastination.com/}. }.

In February 2002, Donald Rumsfeld, the then US Secretary of State for Defence, stated at a Defence Department briefing\footnote{Donald Rumsfeld (February 12, 2002). United States Secretary of Defense. \texttt{https://en.wikipedia.org/wiki/There\_are\_known\_knowns}.}:

\begin{quote}
There are known knowns; there are things we know we know.
We also know there are known unknowns; that is to say we know there are some things we do not know.
But there are also unknown unknowns — there are things we do not know we don't know.
\end{quote}

This thesis had instantly became famous, and when  I first read it, it instantly brought me childhood memories.  I realised that I developed, before Rumsfeld, and used, from a rather early age,   the classification of ``known unknowns'' -- in a context quite different from the one  where Rumsfeld, much later, and independently from me, introduced his famous taxonomy.

When I was about 11 or 12, my parents were almost simultaneously, and suddenly, appointed to posts which, in Western terms, could be described as the  juvenile court judge for our district  (my father) and the child protection officer in the local education authority (my mother) -- this happened because local authorities had to clean up the mess after an especially disgusting scandal of mass sexual abuse of children, and needed two safe pairs of hands.  Inevitably, my parents had to work in close contact, and also inevitably, their decisions were frequently made at the dinner table in our home. And issues were exceptionally private and  sensitive: crimes and criminal prosecution of children and  teenagers, divorces, access rights to, and custody of, children,  adoptions, paternity claims, etc. 

From day one my father told me that I must not only never ever  tell anyone anything that I heard, but  do not even show any sign that I might have heard something, do not even recognise the names of people involved if I heard  them outside of the home. On the other hand, we lived in a small place in deep province, and I inevitably could legitimately know something about someone -- from my school, from our relatives and neighbours,  etc. Moreover, in some cases it was impossible not to know.   

I can proudly say that I had never allowed any slip of tongue or gave any other indication that I was aware of something that I was not supposed to know. To achieve that, I had to learn to classify everything that I\emph{ knew}  about local events and local people in four roughly Rumsfeldian categories (I am using Rumsfeld's   terms  now because I had never had my own names for categories since I had never talked about my classification with anyone). I was quite amused when I first learned about Rumsfeld's taxonomy.

\begin{description}
\item[\textsc{Known Knowns}:]  Stuff that I had legitimate reason to know because everyone knew it, and I knew people involved,  and could freely talk about -- because everyone talked.

\item[\textsc{Unknown Knowns:}]  Stuff that I  was not supposed to know, but about people who I legitimately knew, so I could respond  to someone talking about it:  ``Really? And  how could this happen with them?'' --  although I knew what, when,  how, and why this had happened.

\item[\textsc{Known Unknowns:}]  Stuff that I could legitimately  know (say, some incident everyone was talking about) but was not supposed to know the names of people involved.  When hearing the  story with names mentioned, I had to show surprise, or horror, or some other appropriate emotion and ask something like `` So, you said this was Ivanov? Which one --  the one from the 5th form or his big brother?''
\item[\textsc{Unknown Unknowns:}] I was not supposed to know neither the incident nor people involved.  I had to behave as if I had no vaguest idea what it was and about whom. 
\end{description}

To make things even harder for me, I had to make classification decisions instantly, on the hoof.  

Not long ago, my friend and co-author  \c{S}\"{u}kr\"{u} Yal\c{c}\i nkaya asked me about my parental family. I told him this story, and he said
\begin{quote}
``Aha! Now I know why you always know where you have learned something from, and when.'' 
\end{quote}

It was an exaggeration on his part -- I do not always know, but, indeed, know frequently enough.

\section{December 2020.  Beware of dreams}

In December, lockdown has been  replaced by  Tier 3 restrictions which does  not make much difference. I continue remaining in  state of Nirvana, doing mathematics. My method: to nap several times in a day  and wait for ideas to emerge at awakening -- continues to work. Using results of \cite{borovik-finite} and ideas from my old paper \cite{borovik-periodic-linear} on classification of periodic linear groups I proved a long  standing conjecture \cite[Conjecture 12]{borovik-nesin}. My triumph would be better deserved  if it was someone else's conjecture; unfortunately, it was my own, I formulated it in about 1993.

Perhaps a word of warning is needed: any mathematics produced by a method so irrational as listening to the subconscious  must be immediately, thoroughly, and rigorously checked.

Trust your intuition -- but remember that

\begin{quote}
\emph{verification is the highest form of trust}\footnote{My Russian readers will of course recognise this proverbial saying: \rus{\emph{Проверка -- высшая степень доверия}}.}.
\end{quote}
Verification trains intuition, gives feedback to intuition, to the subconscious. If not properly supported by proofs and calculations, intuition wilts and dies. This simile in a more developed form can be found in my paper \cite[Section 6.1]{borovik-makers-users}, where the training of the subconscious is compared with a training of a dog.

But I would not trust dreams. Perhaps the best authority on that point is Leo Tolstoy. This is from  \emph{War and Peace}, an episode with Pierre Bezukhov' dream\footnote{\emph{War and Peace}, Book Eleven: 1812, Chapter IX.}:

\begin{quotation}


\noindent
``It is dawn,'' thought Pierre. ``But that's not what I want. I want to hear and understand my benefactor's words.'' Again he covered himself up with his cloak, but now neither the lodge nor his benefactor was there. There were only thoughts clearly expressed in words, thoughts that someone was uttering or that he himself was formulating.

Afterwards when he recalled those thoughts Pierre was convinced that someone outside himself had spoken them, though the impressions of that day had evoked them. He had never, it seemed to him, been able to think and express his thoughts like that when awake.

``To endure war is the most difficult subordination of man's freedom to the law of God,'' the voice had said.

``Simplicity is submission to the will of God; you cannot escape from Him. And they are simple. They do not talk, but act. The spoken word is silver but the unspoken is golden. Man can be master of nothing while he fears death, but he who does not fear it possesses all. If there were no suffering, man would not know his limitations, would not know himself. The hardest thing (Pierre went on thinking, or hearing, in his dream) is to be able in your soul to unite the meaning of all. To unite all?'' he asked himself.

 ``No, not to unite. Thoughts cannot be united, but to harness all these thoughts together is what we need! Yes, one must harness them, must harness them!'' he repeated to himself with inward rapture, feeling that these words and they alone expressed what he wanted to say and solved the question that tormented him.

``Yes, one must harness, it is time to harness.''

``Time to harness, time to harness, your excellency! Your excellency!.. some voice was repeating. “We must harness, it is time to harness\dots.''

It was the voice of the groom, trying to wake him. The sun shone straight into Pierre's face. He glanced at the dirty innyard in the middle of which soldiers were watering their lean horses at the pump while carts were passing out of the gate. Pierre turned away with repugnance, and closing his eyes quickly fell back on the carriage seat.

``No, I don't want that, I don't want to see and understand that. I want to understand what was revealing itself to me in my dream. One second more and I should have understood it all! But what am I to do? Harness, but how can I harness everything?'' and Pierre felt with horror that the meaning of all he had seen and thought in the dream had been destroyed.
\end{quotation}

In this translation, the word `\emph{harness}' is used as a translation of two different Russian words of the original: \textbf{sopryagat'} (this what Pierre is heard in his dream) and \textbf{zapryagat}' (this what the groom was actually saying), almost homonyms;  the first is translated as `\emph{conjugate}' and means `\emph{connect}', `\emph{match}'; it is frequently applied to abstract concepts or entities -- and is widely used in mathematics, especially in group theory, my field of study; the second one -- `\emph{harness a horse}'.  `\emph{Conjugate}' could appear in a dream of a mathematician, and mathematicians, I think, should pay attention to  Tolstoy's warning.

And the last but not the least: to retain sanity and the necessary level of confidence in your work, help of colleagues could be useful and even  necessary -- someone else should take a close look at least at the key arguments in your paper. I am lucky that even in lockdown  I can rely on this help from my distant friends.

\section{Boxing Day 2020. Epilogue for the Russian reader}

\begin{flushright}
\emph{The gods have imposed upon my writing the yoke of\\ a foreign tongue that was not
sung at my cradle.}\\ Hermann Weyl
\end{flushright}

A friend who prefers to remain anonymous brought to my attention a wonderful but mostly forgotten Russian poet: Alexander Kochetkov\footnote{\rus{Александр Сергеевич Кочетков} (1900 -- 1953),\\  https://ru.wikipedia.org/wiki/\texttt{\rus{Кочетков\_Александр\_Сергеевич}}.}. His short poem below is about memory and consciousness, it is written in the best tradition of Russian philosophical lyric. After some hesitation and under pressure from my friends I produced a rather clumsy translation of  the subtle and delicate poem -- with my apology, I include it here.

\begin{center}
\parbox[t]{4in}{%
\rus{
\begin{quote}
Мгновенья нет, есть память. Слух полночный\\   Сквозь вздох крови и благовест цветочный \\  Вдруг различит тоскливый некий звук\\   Невидимых орбит (так майский жук\\   Поет под яблоней).  Душа людская,\\   Каким поющим воплем истекая,\\   В какую бездыханность темноты\\   На крыльях памяти несешься ты?..\\
\hspace*{1.55in}  Сергей Кочетков
\end{quote}
}

\smallskip
\begin{quote}
There is  no "now", but  memory exists.\\  The ear of midnight\\ 
Through murmuring of blood and ringing bells of bloom\\
Will suddenly discern a certain wistful sound\\
Of unseen orbits (akin a May-bug's song\\
under an apple tree). A human soul,\\
What is the singing  yell that you are  haemorrhaging?\\
Towards what  kind of suffocating  dark\\
On wings of memory you rashly fly?..\\
\hspace*{1.55in}  Sergey Kochetkov
\end{quote}
}
\end{center}

What we perceive as the stream of our consciousness is memory, rationalization, and ex post facto justification, of decisions and actions already taken, and moreover, started by our subconscious. I heard that neurophysiologists are already making precise measurements in what order and at what intervals different centres in the brain are triggered. We (in the sense of  ``conscious us'') do not make decisions, we justify them retrospectively.

Mathematicians, as it happens with people with a not very clear conscience, have a very sophisticated method of self-justification -- construction of a proof. But even this is done in the same dishonest manner: when we exclaim something like: ``So this should follow from the lemma about the three commutators!'' -- our subconscious has already known this for a long time and only sighs when hears our joyful cry.

Why is it that we are looking for proof with such a morbid obsession? What is our crime that we are hiding from ourselves?

\section*{Acknowledgement}

This text exists because Elizabeth Loew  kindly invited me to contribute a chapter to the book \emph{Math in the Time of Corona}; I decided, out of curiosity, to start a new research project and document its development. This happened to be a success and resulted in the preprint \cite{borovik-finite} and a shorter version of this text \cite{borovik-lockdown}. I decided to expand it after I got more results. It was quite a unusual adventure.

Nothing would be written without help and support from my wife Anna who shares with me life and lockdown.

My work described here is a small fragment of a joint  project with Ay\c{s}e Berkman.\footnote{We are removing `sharpness' from the assumptions of the main result of our paper \cite{berkman-borovik}.} Previously, we did most of our joint work in the Nesin Mathematics Village in \c{S}irince, Turkey,\footnote{https://nesinkoyleri.org/en/nesin-villages/.} and we were planning to meet there during the Easter break and then again in summer\dots\  I also badly missed the planned in advance intensive work sessions with two my other co-authors, Adrien Deloro and \c{S}\"{u}kr\"{u} Yal\c{c}\i nkaya, also supposed to take place in the Village. This is the price I have to pay for the joys of lockdown.

I mentioned here some ideas from the seminal paper by Hall and Higman \cite{hall-higman}. I learned them  in about  1975--76 in the study group which was run by Viktor D. Mazurov  for the benefit of just three undergraduate students:
Elena G. Shapiro, n\'{e}e Bryukhanova,  Evgeny  I.  Khukhro,  and myself -- and, 45 years later, the lockdown gave me a chance to use them, first time in my life.

Special thanks go to  \c{S}\"{u}kr\"{u} Yal\c{c}\i nkaya for explaining to me my past,  to Hayriye B\"{u}\c{s}ra Solak for her comment on Louise Gl\"{u}ck, to my anonymous friend for recommending me the poetry of Alexander Kochetkov, to Ya\v{g}mur Denizhan for persuading me to recklessly translate Kochetkov's poem and to Gregory Cherlin and Adrien Deloro for many enlightening comments and conversations and for expert help at a few critical junctions of my lockdown work.


\section*{Disclaimer} The author writes in his personal capacity and the views expressed do not
necessarily represent the position of his (former) employer or any other person, corporation, organisation, or institution.

\bigskip

\end{document}